\documentclass[12pt]{amsart}

\usepackage{tgpagella}
\usepackage{euler}
\usepackage[T1]{fontenc}
\usepackage{amsmath, amssymb, mathabx}
\usepackage[hidelinks]{hyperref}
\usepackage[english]{babel}
\usepackage{mathrsfs}
\usepackage{eucal}
\usepackage[all]{xy}
\usepackage{tikz}
\usepackage{enumitem}
\newtheorem{thm}{Theorem}[section]
\newtheorem*{thm*}{Theorem}
\newtheorem{lem}[thm]{Lemma}
\newtheorem{fact}[thm]{Fact}

\newtheorem{prop}[thm]{Proposition}
\newtheorem*{prop*}{Proposition}
\newtheorem{conj}[thm]{Conjecture}
\newtheorem{cor}[thm]{Corollary}
\newtheorem*{cor*}{Corollary}

\theoremstyle{definition}
\newtheorem{defn}[thm]{Definition}
\newtheorem*{defn*}{Definition}

\newtheorem{remarks}[thm]{Remarks}
\newtheorem{question}[thm]{Question}

\newtheorem{ex}[thm]{Exercise}

\newtheorem*{question*}{Question}
\newtheorem*{Pquestion*}{Popa's question}

\newtheorem*{conv*}{Convention}

\newcommand{\dminus}{ 
\buildrel\textstyle\ .\over{\hbox{ 
\vrule height3pt depth0pt width0pt}{\smash-} 
}}

\def\bb{\mathbb}

\def\bb{\mathbb}

\def\cal{\mathcal}

\def\u{\mathsf 1}

\newcommand{\cstar}{$\mathrm{C}^*$}

\def \tp{\operatorname{tp}}

\makeatletter

\def\dotminussym#1#2{%
  \setbox0=\hbox{$\m@th#1-$}%
  \kern.5\wd0%
  \hbox to 0pt{\hss\hbox{$\m@th#1-$}\hss}%
  \raise.6\ht0\hbox to 0pt{\hss$\m@th#1.$\hss}%
  \kern.5\wd0}

\DeclareMathOperator{\tr}{tr}

\newcommand\cL{{\cal L}}

\newcommand\bN{{\mathbb N}}

\newcommand\bC{{\mathbb C}}

\def \Th{\operatorname{Th}}
\def \R{\mathcal R}
\def \u{\mathcal U}

\allowdisplaybreaks
\makeatletter
\def\l@subsection{\@tocline{2}{0pt}{2.5pc}{5pc}{}}
\def\l@subsubsection{\@tocline{2}{0pt}{5pc}{7.5pc}{}}
\makeatother

\begin{document}


\title{A survey on the model theory of tracial von Neumann algebras}

\author{Isaac Goldbring and Bradd Hart}
\address{Department of Mathematics\\University of California, Irvine, 340 Rowland Hall (Bldg.\# 400),
Irvine, CA 92697-3875}
\email{isaac@math.uci.edu}
\urladdr{http://www.math.uci.edu/~isaac}
\thanks{Goldbring was partially supported by NSF grant DMS-2054477. Hart was supported by the NSERC}

\address{Department of Mathematics and Statistics, McMaster University, 1280 Main St., Hamilton ON, Canada L8S 4K1}
\email{hartb@mcmaster.ca}
\urladdr{http://ms.mcmaster.ca/~bradd/}

\maketitle

\begin{abstract}
We survey the developments in the model theory of tracial von Neumann algebras that have taken place in the last fifteen years.  We discuss the appropriate first-order language for axiomatizing this class as well as the subclass of II$_1$ factors.  We discuss how model-theoretic ideas were used to settle a variety of questions around isomorphism of ultrapowers of tracial von Neumann algebras with respect to different ultrafilters before moving on to more model-theoretic concerns, such as theories of II$_1$ factors and existentially closed II$_1$ factors.  We conclude with two recent applications of model-theoretic ideas to questions around relative commutants.
\end{abstract}

\section{Introduction}

The tracial ultraproduct construction is an integral tool in the modern study of tracial von Neumann algebras.  Implicitly present in the work of Sakai \cite{sakai}, this construction was of fundamental importance in both the work of McDuff \cite{mcduff} on central sequence algebras and in Connes' landmark result that injective II$_1$ factors are hyperfinite \cite{connes}.   To a model theorist, the presence of an ultraproduct construction makes it natural to try to view tracial von Neumann algebras as structures in an appropriate logic. Discussions around this idea were initiated by Ben Yaacov, Henson, Junge,
and Raynaud during a workshop at the American Institute of Mathematics in 2006 \cite{aimworkshop}.  Farah, Sherman, and the second named author took up this line of inquiry motivated by a question of Popa about isomorphisms of matrix ultraproducts (see Section 3 below).  In a series of three papers \cite{mtoa1,mtoa2,mtoa3}, they introduced a natural language for which the class of tracial von Neumann algebras (as well as the subclasses of tracial factors and II$_1$ factors) become axiomatizable, used model-theoretic ideas to settle a variety of questions about isomorphism of tracial ultraproducts, and initiated the study of elementary equivalence of II$_1$ factors (see Section 4 for a summary of the results proven there).

After these papers, the model theory of operator algebras became a vibrant area of research.  The current authors and Thomas Sinclair \cite{GHS} proved that the class of tracial von Neumann algebras does not admit a model companion.  Later, Farah, Sherman, and the current authors studied the class of existentially closed II$_1$ factors in greater detail \cite{ecfactors}; this work was later complemented by the first author's work on model-theoretic forcing \cite{enforceable} (see also his article in this volume).  These topics are discussed in greater length in Sections 4 and 5.

The general impression is that tracial von Neumann algebras lie on the ``wild'' side of the model-theoretic spectrum.  Indeed, in Section 3 it will be shown that they are almost always unstable and these arguments can be elaborated upon to establish other kinds of unclassifiable Shelah-style model-theoretic behaviour.  On the other hand, model-theoretic techniques can be used to \emph{prove} facts about tracial von Neumann algebras that purely operator-algebraic techniques have been unable to prove.  The last section of this paper exhibits a number of examples of this latter line of research.

We assume that the reader is familiar with the requisite material on tracial von Neumann algebras presented in Adrian Ioana's article in this volume.

\section{The language of tracial von Neumann algebras}

In this section, we describe the language $\cL_{tr}$ for tracial von Neumann algebras together with its intended interpretation. For each $n \in \bN$, there will be a sort $B_n$ together with a symbol $d_n$ for a metric on $B_n$.  If $(M,\tau)$ is a tracial von Neumann algebra, then the intended interpretation of $B_n$ will be the set of all $x \in M$ with $\|x\| \leq n$ (the closed operator norm ball of radius $n$) and for $x,y \in B_n(M)$, the intended interpretation of $d_n(x,y)$ will be $\|x-y\|_2=\sqrt{\tr((x-y)^*(x-y)})$; based on this interpretation, the bound on $d_n$ provided by $\cL_{tr}$ is $2n$.  As we will see, the operator norm will not be a symbol in our language and so there will be work to be done in order to guarantee that the sorts $B_n$ have the intended interpretation in all models of the theory we are describing.  

Regarding the function symbols in $\cL_{tr}$, for every $n \in \bN$, we will have: 
\begin{enumerate}
    \item binary function symbols $+_n$ and $-_n$ with domain $B_n^2$ and range $B_{2n}$.  The intended interpretation is simply addition and subtraction restricted to the operator norm ball of radius $n$.  The modulus of uniform continuity in both cases is the identity function.
    \item a binary function $\cdot_n$ with domain $B_n^2$ and range $B_{n^2}$.  The intended interpretation is the restriction of multiplication to the operator norm ball of radius $n$.  The modulus of uniform continuity of this function symbol is $\epsilon \mapsto \epsilon/n$.
    \item two constant symbols $0_n$ and $1_n$ which lie in the sort $B_n$.  The intended interpretation of these symbols are the elements $0$ and $1$.
    \item for every $\lambda \in \bC$, there is a unary function symbol $\lambda_n$ whose domain is $B_n$ and range is $B_{mn}$, where $m = \lceil{|\lambda|}\rceil$.  The intended interpretation is scalar multiplication by $\lambda$ restricted to the operator norm ball of radius $n$.  The modulus of uniform continuity is $\epsilon \mapsto \epsilon/|\lambda|$.
    \item a unary function symbol $\ast_n$ with domain and range $B_n$.  The intended interpretation is the restriction of the adjoint to the operator norm ball of radius $n$.  The modulus of uniform continuity is the identity function.
    \item for every $m > n$, we include unary function symbols $i_{n,m}$ with domain $B_n$ and range $B_m$.  The intended interpretation is the inclusion map between the balls of the given radii.  The modulus of continuity is the identity map.
\end{enumerate} 

For each $n\in \bb N$, there is one relation symbol (besides the metric symbol) $\tau_n$ whose intended interpretation is the restriction of the trace to the operator norm unit ball of radius $n$.  Its domain is $B_n$ and the range is $\bb D_{n}$, the complex ball of radius $n$.  See the discussion of complex-valued predicates in \cite{}.

Given a tracial von Neumann algebra $M$, we view it is an $\cL_{tr}$-structure as above (interpreting each symbol with its intended meaning); we denote the corresponding $\cL_{tr}$-structure by $D(M)$ and refer to it as the \textbf{dissection of $M$}.

\begin{remarks}

\

\begin{enumerate}
    \item We have introduced the language of tracial von Neumann algebras formally with a subscript for every symbol in a given sort.  In practice, these subscripts will not be used.

 \item The presentation of tracial von Neumann algebras in a continuous logic was superficially different in \cite{mtoa2}.  There, sorts were referred to as ``domains of quantification'' and one did not need connecting maps as everything lay in one common universe. The multi-sorted presentation is more in keeping with other presentations in continuous logic.

\item It is possible to give a presentation of tracial von Neumann algebras with a single sort; for instance, one could use only the operator norm ball of radius 1.  The language would have to change somewhat as, for example, the operator norm unit ball is not closed under addition.

\item It is tempting to consider whether we could present tracial von Neumann algebras in unbounded continuous logic say as presented in \cite{unbdd}.  Although the natural choice of gauge would appear to be that given by the operator norm, since the operator norm is not uniformly continuous with respect to the 2-norm in the class of tracial von Neumann algebras, this na\"ive approach does not appear to work.
\end{enumerate}
\end{remarks}

We now list a set of $\cL_{tr}$-conditions, denoted $T_{tr}$, which we will show axiomatize (dissections of) tracial von Neumann algebras.  We will write these axioms informally and leave it to the reader to write them down as official continuous logic sentences.
\begin{enumerate}
    \item We have equations which tell us that we are dealing with a complex $\ast$-algebra. Because the language is sorted, we need to have an equation for every sort.
    \item There are axioms that express that $\tau$ is a trace.
    \item There are axioms that say that the connecting maps $i_{m,n}$ preserve addition, multiplication, the adjoint and the trace.
    \item $d_n(x,y) = \sqrt{\|x - y\|_2}$. These axioms connect the 2-norm coming from the trace with the metric on each sort and also guarantees that the trace is faithful.
    \item\label{opnorm} $\sup_{x \in B_n} \sup_{y \in B_1} (\|xy\|_2 \dminus n\|y\|_2)$.  These axioms say that elements of the sort $B_n$ have operator norm at most $n$ in the standard representation.
\end{enumerate}


\begin{thm}
$T_{tr}$ axiomatizes the class of (dissections of) tracial von Neumann algebras.
\end{thm}

To prove this theorem, we first note that, given a tracial von Neumann algebra $M$, the dissection $D(M)$ of $M$ is easily seen to be a model of $T_{tr}$. 

Suppose, conversely, that we have a model $A$ of the theory $T_{tr}$.  We begin by forming the direct limit $M$ of the sorts $B_n(A)$ for $n \in \bN$ via the embeddings $i_{m,n}$.  Using the interpretation of the function symbols on each sort, we see that $M$ is naturally a complex $*$-algebra.  Furthermore, using the trace on each sort, one can define an inner product $\langle x,y \rangle := \tau(y^*x)$ on $M$.  We let $H$ be the Hilbert space completion of $M$ with respect to this inner product.  $M$ acts naturally on $H$ by left multiplication. By axiom \ref{opnorm}, left multiplication by any element of $B_n(A)$ has operator norm at most $n$. This allows us to faithfully represent $M$ as a *-subalgebra of $B(H)$.  

We will show that $M$ is actually SOT-closed in $B(H)$ (and is thus a von Neumann algebra) and that $D(M)=A$.  Towards this end, we first observe that any element of $M$ that has operator norm at most $n$ actually belongs to $B_n(A)$; this follows from the same ``functional calculus trick'' that was used in the axiomatization of \cstar-algebras presented in the second author's article in this volume.  Note also that, once we show in the next paragraph that $M$ is SOT-closed, this observation also shows that $D(M)=A$.

We now show that $M$ is SOT-closed.  To see this, suppose that $y\in \overline{M}^{SOT}$; we aim to show that $y\in M$.  Suppose, without loss of generality, that $\|y\|\leq 1$.  By the Kaplansky density theorem, there is a sequence $(x_n)_{n\in \bb N}$ from $M$ with each $\|x_n\|\leq 1$ such that $x_n$ SOT-converges to $y$.  By the previous paragraph, we have that each $x_n\in B_1(A)$.  Now note that $(x_n)$ is Cauchy with respect to the metric $d_1$ on $B_1(A)$.  Indeed, this follows from the fact that 
\[
d_1(x_m,x_n)=\|x_m-x_n\|_2=\|x_m(1)-x_n(1)\|_H\to \|y(1)-y(1)\|_H=0.
\]
Since $B_1(A)$ is $d_1$-complete, it follows that $(x_n)_{n\in \bb N}$ $d_1$-converges to some element $z\in B_1(A)$.  Now given any $w\in M$, we have that $\|x_nw-zw\|_H=\|x_nw-zw\|_2\to 0$; since $M$ is dense in $H$, we have that $x_n$ SOT-converges to $z$, whence $y=z\in M$, as desired.

Once concludes the proof of the theorem by noting that the dissection procedure is an equivalence of categories between tracial von Neumann algebras with embeddings as morphisms and the category of models of $T_{tr}$ again with embeddings as morphisms.

\begin{ex}
Suppose that $(M_i)_{i\in I}$ is a family of tracial von Neumann algebras and $\u$ is an ultrafilter on $I$.  Prove that
$$D\left(\prod_\u M_i\right)=\prod_\u D(M_i),$$ where the ultraproduct on the left-hand side of the equation is the tracial ultraproduct while the ultraproduct on the right-hand side of the equation is the model-theoretic ultraproduct of $\cL_{tr}$-structures.
\end{ex} 

It is known that both the classes of tracial factors and the subclass of II$_1$ factors are closed under (tracial) ultraproducts and ultraroots, and thus form axiomatizable classes.  We can actually write down concrete axioms for these classes.  First, we need the following:

\begin{ex}
Show that the following are definable sets relative to the theory $T_{tr}$:
\begin{itemize}
    \item The set $U$ of unitaries.
    \item The set $P$ of projections.
    \item The set $P_2$ of pairs of projections of the same trace.
\end{itemize}
\end{ex}

The original axiomatization for the class of factors appearing in \cite[Proposition 3.4(1)]{mtoa2} used a Dixmier averaging type argument.  Here, we will instead use the following characterization of factors amongst tracial von Neumann algebras:

\begin{ex}\label{factorex}
Suppose that $M$ is a tracial von Neumann algebra.  Then $M$ is a factor if and only if:  for any $(p,q)\in P_2(M)$, there is $u\in U(M)$ such that $upu^*=q$.
\end{ex}

Consider now the $\cL_{tr}$-sentence $\sigma_{factor}$ given by
$$\sup_{(p,q)\in P_2}\inf_{u\in U}d(upu^*,q).$$  If $M$ is a tracial factor, then $\sigma_{factor}^M=0$ by Exercise \ref{factorex}.  Conversely, if $\sigma_{factor}^M=0$, then given any two projections in $M$ of the same trace, we can a priori only conclude that these projections are approximately unitarily conjugate, meaning that they can be unitarily conjugated arbitrarily close to each other (in 2-norm).  However, by $\aleph_1$-saturation, we see that any two projections in $M^\u$ of the same trace are unitarily conjugate, whence $M^\u$ is a factor by Exercise \ref{factorex}, and, consequently $M$ is also a factor (for otherwise any nontrivial element of the center of $M$ would also be a nontrivial element of the center of $M^\u$).  As a result, we see that $T_{factor}:=T_{tr}\cup\{\sigma_{factor}=0\}$ axiomatizes the class of tracial factors.

Now let $\sigma_{II_1}$ denote the $\cL_{tr}$-sentence $\inf_{p\in P}|\tau(p)-\frac{1}{\pi}|$ and consider the $\cL_{tr}$-theory $T_{II_1}:=T_{factor}\cup\{\sigma_{II_1}=0\}$.  If $M$ is a II$_1$ factor, then $M$ has a projection of trace $\frac{1}{\pi}$, whence $M\models T_{II_1}$.  Conversely, if $M\models T_{II_1}$, then $M$ is a tracial factor that has projections of trace arbitrarily close to $\frac{1}{\pi}$; this precludes $M$ from being type I$_n$ for any $n\geq 1$, whence $M$ is a type II$_1$ factor.  (Of course, nothing is special about $\frac{1}{\pi}$ except that it is irrational.)



\section{Instability and nonisomorphic ultrapowers}

One of the first applications of model theoretic ideas to tracial von Neumann algebras was in addressing the question of the dependence on the ultrafilter of an ultrapower of a given separable tracial von Neumann algebra; this was a question first raised by McDuff in \cite{mcduff}.  To a model theorist, the question of isomorphic ultrapowers is tied up with (model-theortic) \textbf{stability}. We give a precise definition of stability below (Definition \ref{unstabledef}), but the key point is: if $M$ is a separable metric structure in a separable language, then $M$ is stable if and only if all ultrapowers of $M$ with respect to nonprincipal ultrafilters on $\bb N$ are saturated (see \cite[Theorem 5.6]{mtoa2}).  Since any two elementarily equivalent saturated structures are isomorphic (by an easy back-and-forth argument), it follows that all such ultrapowers of a stable structure are isomorphic.

That being said, as we shall soon see, most tracial von Neumann algebras (and, in particular, all II$_1$ factors) are \emph{unstable}.   In the setting of the previous paragraph, when $M$ is unstable, one only knows that the ultrapowers $M^\u$ are $\aleph_1$-saturated and of density character $2^{\aleph_0}$.  Consequently, in models of set theory where the continuum hypothesis holds, such ultrapowers are in fact saturated, the above back-and-forth can still be carried out, and once again one sees that all such ultrapowers are isomorphic.  However, in models of set theory where the continuum hypothesis fails, some ultrapowers will necessarily not be saturated and thus it is possible that there can exist nonisomorphic such ultrapowers.  This possibility is in fact always realized, as we will soon see.

The question of whether a given tracial von Neumann algebra is stable is a function of the type decomposition of the tracial von Neumann algebra.  We recall that from Theorem 3.23 in Ioana's article in this volume that any tracial von Neumann algebra $(M,\tau)$ decomposes as a direct sum $(M,\tau)=(M_1,\tau_1)\oplus (M_2,\tau_2)$, where $M_1$ is a type I tracial von Neumann algebra and $M_2$ is a type II tracial von Neumann algebra.  Moreover, for each $i=1,2$, we can write $M_i=Mz_i$ and $\tau_i=\frac{1}{\tau(z_i)}\tau$, where $z_1,z_2$ are central projections in $M$ such that $z_1+z_2=1$.  It is straightforward to check that, for any ultrafilter $\u$, one has $(M,\tau)^\u=(M_1,\tau_1)^\u\oplus (M_2,\tau_2)^\u$ (see \cite[Lemma 4.1]{mtoa1}).

Let us focus on the type I summand first.  Towards that end, let us suppose that our original $(M,\tau)$ is itself of type I.  In that case, by the discussion in Section 3.5 of Ioana's article in this volume, $M$ further decomposes as $M\cong \bigoplus_n M_n(A_n)$ with each $A_n$ an abelian tracial von Neumann algebra and where $M_n(A_n)$ is equipped with its obvious trace.  Once again appealing to \cite[Lemma 4.1]{mtoa1}, one has $M^\u\cong \bigoplus_n (M_n(A_n))^\u$.  It is immediate to check that $(M_n(A_n))^\u \cong M_n(A_n^\u)$.  This calculation thus reduces the type I situation to the following:

\begin{prop}\label{abeliancase}
Suppose that $(A,\tau)$ is a separable abelian tracial von Neumann algebra.  Then all ultrapowers of $(A,\tau)$ with respect to nonprincipal ultrafilters on $\bb N$ are isomorphic.
\end{prop}

To see this, note that the equivalence of categories between abelian tracial von Neumann algebras and probability algebras is ``ultrapower preserving,'' meaning that it suffices to show that all nonprincipal ultrapowers of a given separable probability algebra $X$ are isomorphic.  One could prove this result using the discussion above together with the fact that any probability algebra is stable (see Berenstein and Henson's article in this volume for a proof of this fact).  However, it will behoove us to give a more direct argument that works in this special case that will also prove useful in our discussion of relative commutants below.

Since the set of atoms do not grow in ultrapowers, one may further assume that $X$ is atomless.  As discussed in Berenstein and Henson's article in this volume, it suffices to show that any nonprincipal ultrapower $X^\u$ of $X$ is Maraham homogeneous of density $2^{\aleph_0}$.

We first show that the top element $1$ has density $2^{\aleph_0}$.  Since it clearly has density at most $2^{\aleph_0}$, it suffices to show that it has density at least $2^{\aleph_0}$.  Towards this end, note that one has a countable family of independent events in $X$, each of measure $\frac{1}{2}$.  It will be beneficial to enumerate these elements by $a_s$ as $s$ ranges over finite subsets of $\bb N$.  Now, given $f:\bb N\to \bb N$, we can consider the element $a_f=(a_{f\upharpoonright n})_\u$ of $X^\u$.  If $f,g:\bb N\to \bb N$ are distinct, then $\mu(a_f\cap a_g)=\frac{1}{4}$, whence $d(a_f,a_g)\geq \frac{1}{4}$ and the density of $1$ is at least $2^{\aleph_0}$, as desired.

To obtain that $X^\u$ is itself homogeneous of density $2^{\aleph_0}$, we have to show that $b$ has denstiy $2^{\aleph_0}$ for any nonzero $b\in X^\u$.  To see this, write $b=(b_n)_\u$ where we may assume each $\mu(b_n)=\mu(b)>0$.  By quantifier-elimination and the strong $\omega$-homogeneity of the theory of atomless probability algebras, there are automorphisms $\sigma_n$ of $X$ for which $\sigma(b_n)=b_0$, yielding an automorphism $\sigma$ of $X^\u$ for which $\sigma(b)=b_0$, where we identify $b_0$ with its diagonal image in $X^\u$.  Now we may repeat the argument from the previous paragraph, this time taking a countable family of independent subsets of $b_0$.

By the discussion preceding Proposition \ref{abeliancase}, we immediately obtain:

\begin{cor}
Suppose that $(M,\tau)$ is a type I tracial von Neumann algebra.  Then all ultrapowers of $(M,\tau)$ with respect to nonprincipal ultrafilters on $\bb N$ are isomorphic.
\end{cor}

So what if $(M,\tau)$ is not of type I?  This is when things get interesting and (in)stability begins to creep in.  We begin with the following simple calculation, which is \cite[Lemma 3.2(1)]{mtoa1}.  In the remainder of this discussion, set $\varphi(x_1,y_1,x_2,y_2)=\|[x_1,y_2]\|_2$, a quantifier-free formula with variables ranging over the unit ball.  

\begin{lem}\label{unstablelemma}
In $M_{2^n}(\bb C)$, there are pairs of contractions $(a_i,b_i)$, $i=1,\ldots,n$, for which $\varphi(a_i,b_i,a_j,b_j)=0$ if $i\leq j$ while $\varphi(a_i,b_i,a_j,b_j)=2$ if $i>j$.
\end{lem}

To see that this lemma holds, write $M_{2^n}(\bb C)=\bigotimes_{i=1}^n M_2(\bb C)$ and set $$a_i=\bigotimes_{k=1}^ix\otimes \bigotimes_{k=i+1}^{n}1, \ b_i=\bigotimes_{k=1}^i1\otimes y\otimes \bigotimes_{k=i+2}^{n}1,$$ where $x=\left(\begin{matrix}0 & \sqrt{2}\\0 & 0\end{matrix}\right)$ and $y=x^*=\left(\begin{matrix}0 & 0\\ \sqrt{2} & 0\end{matrix}\right)$.

The behavior exhibited in the previous lemma is precisely the definition of instability:

\begin{defn}\label{unstabledef}
A metric structure $M$ is \textbf{unstable} (or \textbf{has the order property}) if there is a formula $\psi(x,y)$ and numbers $r < s$ such that: for every $n\in \bb N$, there are tuples $c_i \in M$ for $i < n$ such that $\psi^M(c_i,c_j) \leq r$ for $i\leq j$ while $\psi(c_i,c_j) \geq s$ for $i> j$.
\end{defn}

Using this terminology, Lemma \ref{unstablelemma} implies:

\begin{cor}\label{unstablevna}
Any tracial von Neumann algebra that is not of type I is unstable.
\end{cor}

The previous corollary is indeed a consequence of Lemma \ref{unstablelemma} and the fact that any tracial von Neumann algebra that is not of type I embeds a copy of $M_{2^n}(\bb C)$ for all $n\geq 1$ (see Exercise 5.4 in Ioana's article).  

Returning to the discussion of nonisomorphic ultrapowers, one now quotes the following model-theoretic fact:

\begin{thm}\label{unstablenoniso}
Suppose that the continuum hypothesis fails and $M$ is an unstable separable metric structure.  Then there are nonprincipal ultrafilters $\u$ and $\cal V$ on $\bb N$ such that $M^\u\not\cong M^{\cal V}$.
\end{thm}

The ideas behind this proof are present in Shelah's treatise \cite{shelah} on classification theory but were elucidated in \cite[Proposition 2.6]{mtoa1} (in the special case of tracial von Neumann algebras and related objects) and then in \cite[Theorem 5.6(2)]{mtoa2} (in the general situation).  The basic idea is as follows.  Let $\prec_\psi$ denote the partial order associated with some unstable formula $\psi$, that is, $a\prec_\psi b$ if and only if $\psi(a,b)\leq r$ and $\psi(b,a)\geq s$ (where $r$ and $s$ are as in Definition \ref{unstabledef}).  For an infinite regular cardinal $\lambda$, we define an \textbf{$(\aleph_0,\lambda)$-gap in $M^\u$} to be a pair of sequences $(c_i)_{i<\omega}$ and $(d_\alpha)_{\alpha<\lambda}$ from $M^\u$ such that:
\begin{itemize}
    \item $c_i\prec_\psi c_j$ for all $i<j<\omega$;
    \item $d_\beta\prec_\psi d_\alpha$ for all $\alpha<\beta<\lambda$;
    \item $c_i\prec_\psi d_\alpha$ for all $i<\omega$ and $\alpha<\lambda$;
    \item there \emph{does not} exist $e\in M^\u$ such that $c_i\prec_\psi e\prec_\psi d_\alpha$ for all $i<\omega$ and $\alpha<\lambda$.
\end{itemize}
Set $\kappa(M^\u)$ to be the minimal cardinal $\lambda$ such that there is an $(\aleph_0,\lambda)$-gap in $M^\u$.  Note that, by instability and $\aleph_1$-saturation, we have that $\kappa(M^\u)\geq \aleph_1$.  Also, let $\kappa(\u)$ denote the coinitiality of the set of infinite elements in the (ordinary, discrete model theory) ultrapower $(\bb N,<)^\u$ of the structure $(\bb N,<)$.   The proof of Theorem \ref{unstablenoniso} relies on the fact that (i) $\kappa(M^\u)=\kappa(\u)$, and (ii) under the negation of the continuum hypothesis, there are nonprincipal ultrafilters $\u$ and $\cal V$ on $\bb N$ such that $\kappa(\u)\not=\kappa(\cal V)$.  In fact, for each regular cardinal $\aleph_1\leq \kappa\leq 2^{\aleph_0}$, there is a nonprincipal ultrafilter $\u$ on $\bb N$ such that $\kappa(\u)=\kappa$.  (Note, in particular, that when $\kappa(\u)<2^{\aleph_0}$, the ultrapower $M^\u$ is not saturated.)  By a result of Solovay, it is consistent that there are continuum many $\kappa$ in between $\aleph_1$ and $2^{\aleph_0}$, whence one would have continuum many such nonisomorphic ultrapowers.  Farah and Shelah \cite{farahshelah} improved this result by showing that one can in fact arrange for $2^{2^{\aleph_0}}$ many nonisomorphic such ultrapowers of a given unstable separable structure (assuming the failure of the continuum hypothesis).  

Proposition \ref{abeliancase}, Corollary \ref{unstablevna}, and Theorem \ref{unstablenoniso} immediately imply:

\begin{cor}
For a given separable tracial von Neumann algebra $(M,\tau)$, we have the following:
\begin{itemize}
    \item If $(M,\tau)$ is type I, then $(M,\tau)$ is stable and all nonprincipal ultrapowers of $(M,\tau)$ are isomorphic.
    \item If $(M,\tau)$ is not type I, then $(M,\tau)$ is unstable.
    \begin{itemize}
    \item If the continuum hypothesis holds, then all nonprincipal ultrapowers of $(M,\tau)$ are isomorphic.
        \item If the continuum hypothesis fails, then there are nonprincipal ultrafilters $\u$ and $\cal V$ on $\bb N$ such that $(M,\tau)^\u\not\cong (M,\tau)^{\cal V}$.
    \end{itemize}
\end{itemize}
\end{cor}

One can use a variant of the above argument to answer a question of Popa:

\begin{thm}\label{popanoniso}
Suppose that the continuum hypothesis fails.  Then there are nonprincipal ultrafilters $\u$ and $\cal V$ on $\bb N$ such that $\prod_\u M_n(\bb C)\not\cong \prod_{\cal V}M_n(\bb C)$.
\end{thm}

To deduce the previous theorem, one first argues that, given any $\epsilon>0$ and $n\in \bb N$, there is $m_0\in \bb N$ such that, for all $m\geq m_0$, one can find pairs of contractions $(a_i,b_i)$, $i=1,\ldots,n$, in $M_m(\bb C)$ such that $\varphi(a_i,b_i,a_j,b_j)\leq \epsilon$ if $i\leq j$ while $\varphi(a_i,b_i,a_j,b_j)\geq 2-\epsilon$ if $i>j$.  The key to doing that is to take $m_0$ sufficiently larger than $2^n$ so that, for $m\geq m_0$, writing $m=k2^n+r$ with $0\leq r<2^n$, the fraction $\frac{r}{m}$ is neglibible.  Now taking $p$ a projection in $M_m(\bb C)$ of trace $\frac{k2^n}{m}$, one has that $pM_m(\bb C)p\cong M_{2^n}(\bb C)\otimes M_k(\bb C)$ and the contractions $(a_i,b_i)\in M_{2^n}(\bb C)$ guaranteed by Lemma \ref{unstablelemma} are the desired elements of $M_m(\bb C)$ under this identification.  This fact is enough to still find an $(\aleph_0,\lambda)$-gap in $\prod_\u M_n(\bb C)$ and the minimal such $\lambda$ is once again equal to $\kappa(\u)$, allowing one to conclude as before.

It is interesting to note that the conclusion of Theorem \ref{popanoniso} might still be true even if the continuum hypothesis is assumed.  Indeed, we do not know if the following question has a positive answer:

\begin{question}
Do we have $\prod_\u M_n(\bb C)\equiv \prod_{\cal V}M_n(\bb C)$ for all nonprincipal ultrafilters $\u$ and $\cal V$ on $\bb N$?  Said differently:  does $\lim_{n\to \infty}\sigma^{M_n(\bb C)}$ exist for all sentences $\sigma$?
\end{question}

The above considerations can also be used to shed light on a question of McDuff \cite{mcduff}, who asked:  given a separable II$_1$ factor $M$, does the isomorphism type of the relative commutant $M'\cap M^\u$ depend on the choice of the nonprincipal ultrafilter $\u$?  McDuff had showed that there were three possibilities for the type of the relative commutant:  (i)  it equals $\bb C$ (in which case $M$ does not have property Gamma, see Theorem 6.19 in Ioana's article);  (ii) it is abelian, atomless, and has density character $2^{\aleph_0}$ (in which case $M$ does have property Gamma but is not McDuff, see Theorem 6.23 in Ioana's article); or (iii) it has type II$_1$ (and is thus McDuff).  We stress that the category that the relative commutant $M'\cap M^\u$ falls into is independent of the choice of $\u$.  In that case, it is clear that the isomorphism type of the relative commutant does not depend on $\u$ when $M$ does not have property Gamma.

In the case that $M$ has property Gamma but is not McDuff, we argue in a manner similar to our proof of Theorem \ref{abeliancase} above.  We first show that $M'\cap M^\u$ has density character $2^{\aleph_0}$.  Once again, the point is to show that it has density character at least $2^{\aleph_0}$.  Since the probability algebra associated to $M'\cap M^\u$ is atomless, we can find countably many commuting projections in $M'\cap M^\u$ of trace $\frac{1}{2}$.  This allows us to find, for finite subsets $s$ of $\bb N$, ``almost'' projections $p_s$ in $M$ of trace ``almost'' $\frac{1}{2}$, which ``almost'' commute with each other and which ``almost'' commute with $M$, and for which products $p_sp_t$ have trace ``almost'' $\frac{1}{4}$ for distinct $s,t$.  This allows us to define, for $f:\bb N\to \bb N$, elements $p_f=(p_{f\upharpoonright n})_\u\in M'\cap M^\u$.  For distinct $f,g:\bb N\to \bb N$, one has that $\tau(p_fp_g)=\frac{1}{4}$, whence the elements $p_f$ represent a discrete subset of $M'\cap M^\u$ of cardinality $2^{\aleph_0}$, as desired.

To finish the argument, it suffices to show that the same result holds below any given nonzero projection $p\in M'\cap M^\u$.  One many write $p=(p_n)_\u$ where, without loss of generality, we may assume that all $p_n$'s have the same trace as $p$.  One then argues as before, finding elements $p_s\in M^\u$ satisfying all of the requirements above as well as the extra requirement that $p_s\leq p_n$ if $s\subseteq \{0,1,\ldots,n-1\}$.  The resulting elements $p_f$ now lie in $M'\cap M^\u$ and are below $p$. 

Finally, assume that $M$ is McDuff.  In this case, one can find a copy of the hyperfinite II$_1$ factor $\cal R$ in $M'\cap M^\u$.  This allows one to witness the instability above in $M'\cap M^\u$.  One can then analogously define and find $(\aleph_0,\lambda)$-gaps in $M'\cap M^\u$ and the minimal length $\lambda$ of such a gap once again coincides with $\kappa(\u)$, allowing one to conclude as above.  The details of this argument can be found in \cite[Proposition 2.9]{mtoa1}.

\section{Theories of II$_1$ factors}

As we know from the previous section, II$_1$ factors are wild from the stability-theoretic point of view.  In fact, \cite[Theorem 4.3]{mtoa3} shows that any II$_1$ factor has the maximum number of separable models:

\begin{thm}[Farah, Hart, and Sherman]\label{manymodels}
For any separable II$_1$ factor $N$, there are continuum many nonisomorphic separable II$_1$ factors with the same theory as $N$.
\end{thm}

The proof of the previous proposition uses a result of Nicoara, Popa, and Sasyk \cite{NPS}, which constructs a family $(M_\alpha)_{\alpha\in 2^\omega}$ of separable II$_1$ factors, each of which embeds in $\R^\u$, and for which at most countably many embed into any given separable II$_1$ factor.  Since $\R$ embeds into $N$, we see that $M_\alpha$ embeds into $N^\u$ for each $\alpha$.  Let $N_\alpha$ be a separable elementary substructure of $N^\u$ containing the image of $M_\alpha$.  Then each $N_\alpha$ has the same theory as $N$ and there are uncountably many nonisomorphic $N_\alpha$'s by the defining property of the family $(M_\alpha)_{\alpha\in 2^\omega}$.

For some II$_1$ factors $N$ , one can find a more ``concrete'' construction of uncountably many nonisomorphic seprable models of $\Th(N)$, as discussed right before \cite[Theorem 4.3]{mtoa3}.  Recall from Section 6.2 of Ioana's article that the fundamental group of a II$_1$ factor $N$ is the set $\cal F(N):=\{t\in \bb R^{>0} \ : \ N\cong N_t\}$, where $N_t$ is the amplification of $N$ by $t$.  $\cal F(N)$ is indeed a subgroup of the multiplicative group $\bb R^{>0}$ (whence the name).  Let $N$ be a II$_1$ factor whose fundamental group is a countable dense subgroup of $\bb R^{>0}$ (such a II$_1$ factor was first constructed by Golodets-Nessonov in \cite{GN}).  Then the index of $\cal F(N)$ in $\bb R^{>0}$ is $2^{\aleph_0}$ and for $s,t\in \bb R^{>0}$ in distinct cosets of $\cal F(N)$, we have that $N_s\not\cong N_t$.  However, as shown in \cite{GHgames}, the \textbf{first-order fundamental group} $\cal F_{fo}(N):=\{t\in \bb R^{>0} \ : \ N\equiv N_t\}$ is a \emph{closed} subgroup of $\bb R^{>0}$ containing $\cal F(N)$, whence, in our case, must be all of $\bb R^{>0}$, that is, each $N_t$ is elementarily equivalent to $N$.  


We now turn our attention towards finding different theories of II$_1$ factors.  Interestingly enough, the history of finding non-elementarily equivalent II$_1$ factors followed the history of finding non-isomorphic separable II$_1$ factors.  The first example of non-isomorphic separable II$_1$ factors was in Murray and von Neumann's paper \cite{MvN}, where they showed that $\R$ has property Gamma while $L(\bb F_2)$ does not.  In \cite[Subsection 3.2.2]{mtoa3}, it was shown that property Gamma is an $\forall\exists$-axiomatizable property of II$_1$ factors, whence $\R$ and $L(\bb F_2)$ are not even elementarily equivalent.

In connection with the previous paragraph, it is somewhat surprising that the following question (first raised in \cite[Question 3.5]{mtoa3}) is still open:

\begin{question}
Is the class of II$_1$ factors that do \emph{not} have property Gamma an elementary class?
\end{question}

We digress for a moment to discuss the general progress in understanding first-order theories of II$_1$ factors without property Gamma.  In general, this line of inquiry is not well understood.  In fact, the first known proof that there are at least two elementary equivalence classes of II$_1$ factors without property Gamma used the fact that there is a counterexample to the Connes Embedding Problem without property Gamma!  More precisely, if $M$ is a II$_1$ factor that does not embed into $\R^\u$, then the free product $M*\R$ is a II$_1$ factor without property Gamma which also does not embed into $\R^\u$, and thus cannot have the same theory (even universal theory) as, say, $L(\bb F_2)$.  More recently, in \cite{cik}, Chifan, Ioana, and Kunnawalkam Elayavalli gave an explicit construction (not reliant on the resolution of the Connes Embedding Problem) of a II$_1$ factor $N$ without property Gamma which is not elementarily equivalent to $L(\bb F_2)$.  There, one can deduce that $N$ must have a different $\forall\exists$-theory from $L(\bb F_2)$; it is unknown whether or not the $N$ constructed there embeds in $\R^\u$.

In light of the previous paragraph, the following question is still open:

\begin{question}\label{nongammaquestion}
Are there II$_1$ factors $M$ and $N$, neither of which have property Gamma, both of which embed into $\R^\u$, for which $M\not\equiv N$?
\end{question}

One might attempt to settle the previous question by considering the following natural model-theoretic variant of the well-known free group factor problem:

\begin{question}
Are $L(\bb F_m)$ and $L(\bb F_n)$ elementarily equivalent for distinct $m,n\geq 2$?
\end{question}

In \cite[Theorem 5.1]{mtoa3}, it is shown that any nonprincipal ultraproduct of matrix algebras does not have property Gamma.  In connection with Question \ref{nongammaquestion}, it is thus interesting to ask:

\begin{question}
Are there $n\geq 2$ and a nonprincipal ultrafilter $\u$ on $\bb N$ for which $L(\bb F_n)\equiv \prod_\u M_n(\bb C)$?
\end{question}

We point out Chifan, Ioana, and Kunnawalkam Elayavalli show that their II$_1$ factor $N$ without property Gamma that is not elementarily equivalent to $L(\bb F_2)$ is also not elementarily equivalent to $\prod_\u M_n(\bb C)$ for any nonprincipal ultrafilter on $\bb N$. 



We now return to the main thread concerning new theories of II$_1$ factors with property Gamma.  The next example of a ``new'' II$_1$ factor used the property of being McDuff.  It is trivial to see that $\R$ is a McDuff II$_1$ factor.  Dixmier and Lance \cite{DL} produced an example of a II$_1$ factor $M$ that had property Gamma (so $M\not\cong L(\bb F_2)$) but was not McDuff (whence $M\not\cong \R$).  In \cite[Proposition 3.9]{mtoa3}, it was shown that being McDuff is an $\forall\exists$-axiomatizable property of separable II$_1$ factors (and these axioms can be used to define McDuff II$_1$ factors of arbitrary density character), whence the Dixmier and Lance factor represents a third elementary equivalence class of II$_1$ factors.

Throughout the years, a handful of new examples of separable II$_1$ factors appeared, albeit at a fairly slow pace.  (It is worth mentioning that one of these new examples, introduced by Zeller-Meier in \cite{ZM}, also represented a new, that is, fourth, elementary equivalence class, as pointed out by the authors in \cite{GHgames}.)  A breakthrough finally occurred in 1969, where in a pair of papers \cite{McD1} and \cite{McD2}, McDuff constructed first countably many and then continuum many nonisomorphic separable II$_1$ factors.  It is worth pointing out that McDuff's construction is quite explicit.  Indeed, she constructs two (very concrete) functors $T_0$ and $T_1$ from the category of countable groups into itself, which can then be composed, yielding functors $T_\alpha$ for every $\alpha\in 2^{<\omega}$.  If $\alpha$ is an initial segment of $\beta$, it is fairly easy to see that $T_\alpha(G)$ is a subgroup of $T_\beta(G)$, whence one can define $T_\alpha(G)$ for any $\alpha\in 2^{\omega}$ by taking the direct limit of the groups $T_{\alpha|n}(G)$ for $n\in \omega$.  McDuff's examples are then defined by $M_\alpha:=L(T_\alpha(\bb F_2))$ for $\alpha\in 2^\omega$.

Amazingly enough, Boutonnet, Chifan, and Ioana \cite{BCI} showed that McDuff's factors are not even elementarily equivalent, whence showing that there are continuum many distinct theories of II$_1$ factors (a fact suspected by those working in the area):

\begin{thm}[Boutonnet, Chifan, and Ioana]\label{BCItheorem}
Let $(M_\alpha)_{\alpha<2^{\omega}}$ denote the family of continuum many nonisomorphic separable II$_1$ factors constructed by McDuff in \cite{McD2}.  Then for any ultrafilters $\u$ and $\cal V$ (on arbitrary index sets) and any distinct $\alpha,\beta\in 2^{\omega}$, we have $M_\alpha^\u\not\cong M_\beta^{\cal V}$.
\end{thm}

Using the Keisler-Shelah theorem for continuous logic, it follows that $M_\alpha\not\equiv M_\beta$ for any distinct $\alpha,\beta\in 2^{\omega}$.  However, to a model theorist, this method of proof for demonstrating continuum many distinct theories of II$_1$ factors is not quite satisfying.  Indeed, it would be preferable to point to a single sentence (or perhaps a small, concrete set of sentences) which distinguishes the theory of $M_\alpha$ from the theory of $M_\beta$.  At first glance, it was not clear at all how to achieve such a goal, for Boutonnet, Chifan, and Ioana produced invariants of the ultrapowers (and not the separable objects themselves) that distinguished the ultrapowers of distinct McDuff factors.

That being said, using Ehrenfeucht-Fra\"isse games and a careful analysis of the proof of Theorem \ref{BCItheorem}, the current authors were at least able to find an upper bound on the quantifier-complexity of sentences that could be used to distinguish distinct McDuff factors:

\begin{thm}[Goldbring and Hart \cite{GHgames}]
Suppose that $\alpha,\beta\in 2^{\omega}$ are distinct sequences and let $k<\omega$ be minimal such that $\alpha(k)\not=\beta(k)$.  Then $M_\alpha\not\equiv_{5k+3}M_\beta$.
\end{thm}

Later, with Towsner \cite{GHT}, using a much more careful analysis, the authors were able to find (essentially) concrete sentences which distinguished the McDuff factors.  It is satisfying to note that the quantifier-complexity of these sentences agreed with the theoretic bounds predicted by the previous theorem.  Implicit in the arguments in \cite{GHgames} and \cite{GHT} is the theory of definability in continuous logic and the connection with definability (and the operator-algebraic notion of spectral gap, see Subsection \ref{FCEP} below) was highlighted and made much more explicit in the first author's article \cite{spgap}.

As pointed out by the authors in \cite{GHgames}, a positive answer to the following question would yield an alternative proof of the existence of continuum many different theories of II$_1$ factors:

\begin{question}
Is there a II$_1$ factor $M$ such that the first-order fundamental group $\cal F_{fo}(M)$ of $M$ is a \emph{proper} subgroup of $\bb R^{>0}$? 
\end{question}

\section{Existentially closed II$_1$ factors and nonexistence of model companions}

A typical model-theoretic inquiry into a natural (elementary) class of structures often begins by trying to prove some form of quantifier elimination or model completness result, that is, by attempting to prove the existence of a model companion or model completion for this class.  As is to be expected from the wild model-theoretic behavior of tracial von Neumann algebras discussed in the previous two sections, it should not be too surprising that the model companion for the class of tracial von Neumann algebras does not exist, as was first shown by the authors and Thomas Sinclair in \cite{GHS}.  The starting point in proving this result is the following:

\begin{thm}[Goldbring, Hart, and Sinclair]\label{noqe}
$\Th(\R)$ does not have quantifier-elimination.
\end{thm}

To prove this result, it suffices to find a model $M$ of $\Th_\forall(\R)$, an automorphism $\alpha\in \operatorname{Aut}(M)$, and an embedding $\iota:M\hookrightarrow \R^\u$ for which $\alpha$ does not extend to an inner automorphism of $\R^\u$, that is, for which there does not exist a unitary $u\in U(\R^\u)$ for which $\iota(\alpha(x))=u\iota(x)u^*$ for all $x\in M$.  Indeed, suppose this can be done and set $N:=M\rtimes_\alpha \bb Z$ to be the crossed product algebra (see Section 3.6 of Ioana's article for details on this construction).  In particular, $N$ is a tracial von Neumann algebra containing $M$ which contains a unitary $v$ that implements $\alpha$, that is, $\alpha(x)=vxv^*$ for all $x\in M$.  It follows that there can be no embedding $\iota':N\hookrightarrow \R^\u$ extending $\iota$ (lest $\iota'(v)\in U(\R^\u)$ contradict the choice of $M$, $\alpha$, and $\iota$), whence, since it can be shown that $N$ is also a model of $\Th_\forall(\R)$, it follows that $\Th(\R)$ does not have quantifier-elimination.  It remains to notice that Brown constructed such objects in \cite{nate}; in fact, one can take $M:=L(\operatorname{SL}_3(\bb Z)*\bb Z)$ and $\alpha$ any automorphism of the form $\operatorname{id}*\theta$ where $\theta\in L(\bb Z)$ is nontrivial.  (The description of $\iota$ is a bit more complicated to describe.)

It is interesting to note that any matrix algebra $M_n(\bb C)$ does in fact admit quantifier-elimination; see \cite[Theorem 2.2]{jing}.  (We thank David Jekel for pointing this out to us, although this was apparently also known to Ilijas Farah and others.)

After the fact, one realizes that the Brown result used in the proof of Theorem \ref{noqe} works equally well for any McDuff II$_1$ factor, whence one sees that $\Th(M)$ does not admit QE for any McDuff II$_1$ factor.  (In \cite{GHS} it was also assumed that $M$ is also a so-called locally universal factor, but this is unnecessary.)  The relevance of this observation comes from the following:

\begin{lem}\label{ecvnamcduff}
Any existentially closed tracial von Neumann algebra is a McDuff II$_1$ factor.
\end{lem}

Recall that if $\cL$ is any continuous language and the $\cL$-structure $M$ is a substructure of the $\cL$-structure $N$, then $M$ is said to be \textbf{existentially closed in $N$} (or \textbf{e.c. in $N$} for short) if:  for any quantifier-free formula $\varphi(x,y)$ and parameters $b\in M$, we have
$$\inf\{\varphi(a,b)^M \ : \ a\in M\}=\inf\{\varphi(a,b)^N \ : \ a\in N\}.$$ Alternatively, $M$ is e.c. in $N$ if and only if there is an embedding $N\hookrightarrow M^\u$ that restricts to the diagonal embedding of $M\hookrightarrow M^\u$.  If $T$ is an $\cL$-theory, we say that $M$ is an \textbf{e.c. model of $T$} if $M$ is a model of $T$ and is e.c. in any extension that is a model of $T$.  Consequently, by an existentially closed tracial von Neumann algebra, we mean an existentially closed model of $T_{tr}$.

The proof of Lemma \ref{ecvnamcduff} follows from (i) the class of McDuff II$_1$ factors is $\forall\exists$-axiomatizable, and (ii) any tracial von Neumann algebra embeds into a McDuff II$_1$ factor.  To see the latter, note that if $M$ is a tracial von Neumann algebra, then $M*L(\bb Z)$ is a II$_1$ factor and then $(M*L(\bb Z))\otimes \R$ is a McDuff II$_1$ factor.

Recall that a continuous theory $T$ \textbf{has a model companion} if the class of e.c. models of $T$ is an elementary class (equivalently, if the class of e.c. models is closed under ultraproducts).  In this case, the model companion is unique (up to logical equivalence) and is the theory axiomatizing the class of e.c. models of $T$.

Armed with all this, one obtains quite easily:

\begin{cor}[Goldbring, Hart, and Sinclair]\label{nomodelcompanion}
The theory of tracial von Neumann algebras does not have a model companion.
\end{cor}

Indeed, since the class of tracial von Neumann algebras has the amalgamation property (as witnessed by the amalgamated free product construction, a generalization of the free product construction described in Section 5.6 of Ioana's article), any model companion would be a \textbf{model completion}, that is, would have quantifier-elimination.  However, by the previous lemma, the models of the model completion would be McDuff II$_1$ factors, which we just argued can never have quantifier-elimination.

An alternative proof of Corollary \ref{nomodelcompanion} was given by the first author in \cite[Corollary 5.19]{spgap} and did not rely on Brown's work (but ultimately simply rests on the the existence of property (T) groups, which is itself part of Brown's work referred to above).

Recall the famous \textbf{Connes Embedding Problem } (CEP) asks whether or not all II$_1$ factors have the same universal theory, namely $\Th_\forall(\R)$.  We now know that the CEP has a negative solution; see the first author's article in this volume for much, much more on the CEP.  That being said, it thus becomes interesting to ask whether or not $\Th_\forall(\R)$ admits a model companion.  This question also has a negative answer, as was first shown by Farah, Sherman, and the authors in \cite{ecfactors} (years before the negative solution to the CEP was announced).   The first thing to note is the following:

\begin{lem}\label{Rec}
$\R$ is an e.c. model of $\Th_\forall(\R)$.
\end{lem}

This lemma follows immediately from the following well-known fact about $\R$:

\begin{fact}\label{folklore}
Any two embeddings $\rho_1,\rho_2:\R\hookrightarrow \R^\u$ are \emph{unitarily conjugate}, that is, there is a unitary $u\in U(\R^\u)$ such that $\rho_1(x)=u\rho_2(x)u^*$ for all $x\in \R$.
\end{fact}

Indeed, a consequence of this fact is that any embedding of $\R$ in its ultrapower, being unitarily conjugate to the digaonal embedding, is elementary; Lemma \ref{Rec} follows immediately from this fact.  The fact that any embedding of $\R$ into its ultrapower is elementary shows that $\R$ is a \textbf{finitely generic} model of $\Th_\forall(\R)$; see \cite{enforceable}.

While Fact \ref{folklore} is quite well-known, it is difficult to find in print, so a word about its proof might be in order.  First, one uses an argument involving Murray-von Neumann equivalence of projections in II$_1$ factors to show that any two embeddings $M_n(\bb C)\hookrightarrow N$ (for $N$ a II$_1$ factor) are unitarily conjugate.  A standard approximation argument then shows that any two embeddings $\rho_1,\rho_2:\R\hookrightarrow N$ (again, for $N$ a II$_1$ factor) are \emph{approximately unitarily conjugate}, meaning that for any finite set $F\subseteq \R$ and any $\epsilon>0$, there is $u\in U(N)$ such that $\|\rho_1(x)-u\rho_2(x)u^*\|_2<\epsilon$ for all $x\in F$.  If $N$ is additionally assumed to be $\aleph_1$-saturated (e.g. when $N$ is an ultraproduct), then this approximate unitary conjugacy can be upgraded to actual unitary conjugacy.  A model-theoretic proof of Fact \ref{folklore} was given by the first author in \cite{resembles}.

Lemma \ref{Rec} might lead one to believe that $\Th(\R)$ is model-complete and is the model companion (but not model completion by Theorem \ref{noqe}) of its universal theory.  As alluded to above, this is in fact not the case:

\begin{thm}[Farah, Goldbring, Hart, and Sherman \cite{ecfactors}]\label{notmodelcomplete}
$\Th(\R)$ is not model-complete.  In other words, $\Th_\forall(\R)$ does not have a model companion.
\end{thm}

We note that one can not simply prove this theorem as in the unrestricted case above as the following question is currently open:

\begin{question}\label{amalgamated}
Does $\Th_\forall(\R)$ have the amalgamation property?  In particular, if $M$ and $N$ are models of $\Th_\forall(\R)$ with a common subalgebra $Q$, is the amalgamated free product $M*_QN$ once again a model of $\Th_\forall(\R)$?
\end{question}

Surprisingly, the state of knowledge of the latter question is that the answer is positive if one assumes that the algebra $Q$ is \emph{amenable} (a very serious restriction).  We mention in passing that a positive answer to Question \ref{amalgamated} yields a positive answer to Question \ref{nongammaquestion} above.  Indeed, the II$_1$ factor without property Gamma that is also not elementarily equivalent to $L(\bb F_2)$ explicitly constructed by Chifan, Ioana, and Kunnawalkam Elayavalli mentioned above does in fact embed into $\R^\u$ provided Question \ref{amalgamated} has a positive answer.

The proof of Theorem \ref{notmodelcomplete} rests on the following fundamental result of Jung, which gives a converse to Fact \ref{folklore}:

\begin{fact}[Jung \cite{jung}]\label{jungfact}
Suppose that $M$ is a finitely generated II$_1$ factor that is a model of $\Th_\forall(\R)$.  Further suppose that any two embeddings $M\hookrightarrow \R^\u$ are unitarily conjugate.  Then $M\cong \R$.
\end{fact}

The proof of Theorem \ref{notmodelcomplete} proceeds by showing that, under the assumption that $\Th(\R)$ is model-complete, any ``nonstandard'' model of $\Th(\R)$, that is, any $\R'\equiv \R$ with $\R'\not\cong \R$ (which exists by Theorem \ref{manymodels}) satisfies the assumption of Fact \ref{jungfact}, yielding the desired contradiction.  To see this, we first note that $\R'$ is finitely generated (in fact, singly generated) since it is a McDuff II$_1$ factor; fix a generator $a$ of $\R'$.  Now fix embeddings $j_1,j_2:\R'\hookrightarrow \R^\u$; we aim to show that $j_1$ and $j_2$ are unitarily conjugate, that is, we wish to find $u\in U(\R^\u)$ such that $j_2(a)=uj_1(a)u^*$.  Since $\Th(\R)$ is model-complete, $j_1$ and $j_2$ are elementary embeddings, whence $\tp^{\R^\u}(j_1(a))=\tp^{\R^\u}(j_2(a))$.  Take an elementary extension $\bar \R$ of $\R^\u$ and an automorphism $\alpha$ of $\bar R$ such that $\alpha(j_1(a))=j_2(a)$.  By considering $N:=\bar \R\rtimes_\alpha \bb Z$, which contains a unitary implementing $\alpha$, the fact that $\R^\u$ is existentially closed (which follows from $\Th(\R)$ being model-complete) and $\aleph_1$-saturated implies the existence of the desired unitary.

An alternative proof of Theorem \ref{notmodelcomplete}, showing directly that $\R^\u$ is not an existentially closed model of $\Th_\forall(\R)$, is given in \cite[Theorem 2.3.4]{AGK}; the proof there relies on deeper results but still is of interest in its own right.

Despite these negative results, understanding properties of e.c. II$_1$ factors has been important in applications of model-theoretic techniques in studying purely operator-algebraic problems, as we will see in the next section.  We are still far from having a great understanding of the class of e.c. II$_1$ factors.  In fact, the following question is still open: 

\begin{question}\label{eceequestion}
Do there exist non-elementarily equivalent e.c. II$_1$ factors?  Do there exist non-elementarily equivalent e.c. models of $\Th_\forall(\R)$?
\end{question}

The analogous question for e.c. groups was famously solved by Macintyre in \cite{macintyre} using model-theoretic forcing; it is our hopes that such techniques might also prove useful in settling the previous question.

Some further properties of e.c. factors can be found in Chifan, Drimbe, and Ioana's article in this volume.

\section{Two recent applications}

In this section we discuss two recent applications of model-theoretic techniques in von Neumann algebra theory.

\subsection{Popa's Factorial Commutant Embedding Problem}\label{FCEP}

Popa formulated the following problem in connection with the CEP:

\begin{conj}[Popa's Factorial Commutant Embedding Problem (FCEP)]\label{popafcep}
Suppose that $N$ is a separable II$_1$ factor that embeds into $\R^\u$.  Further suppose that $N$ has \emph{property (T)}.  Then there is an embedding $i:N\hookrightarrow \R^\u$ such that the relative commutant $i(N)'\cap R^\u$ is a factor. 
\end{conj}

Here, property (T) is a certain rigidity property of a von Neumann algebra (see Section 6.2 in Ioana's article in this volume).  A large class of examples of von Neumann algebras with property (T) are the group von Neumann algebras $L(G)$ corresponding to a discrete property (T) group $G$, for example $G=\operatorname{SL}_3(\bb Z)$.  The FCEP has a positive solution for $N=\operatorname{SL}_3(\bb Z)$ (as Popa himself showed in \cite[Section 1.7]{popa}).  Besides this result, very little progress had been made on this problem.

The FCEP has a nice geometric interpretation and connects with the work of Brown mentioned in the previous section.  Indeed, for any model $M$ of $\Th_\forall(\R)$, Brown considers the set $\bb Hom(M,\R^\u)$ of embeddings of $M$ into $\R^\u$ modulo unitary conjugacy.  He endows this set with a natural ``convex-like structure'' and then shows that the extreme points of this space are precisely those embeddings whose image has factorial commutant.  Thus, the FCEP asks if every such space $\bb Hom(M,\R^\u)$ has extreme points.  

In \cite{fcep}, the first author made some progress on the FCEP by proving the following theorem:

\begin{thm}\label{goldfcep}
There is a \emph{locally universal} II$_1$ factor $M$ such that every separable property (T) factor $N$ admits an embedding $i:N\hookrightarrow M^\u$ for which $i(N)'\cap M^\u$ is a factor.
\end{thm}

Here, the fact that $M$ is locally universal means every II$_1$ factor admits an embedding into an ultrapower of $M$.  By the negative solution of the CEP referred to above, we know that $M$ cannot be $\R$ as in the original version of the FCEP (nor can it even embed in $\R^\u$); nevertheless, as we will discuss momentarily, the proof of Theorem \ref{goldfcep} may indicate how one might settle Conjecture \ref{popafcep} in its entirety.

In what follows, the only particular aspect of property (T) von Neumann algebras we will use is the following:

\begin{defn}
If $N$ is a subalgebra of $M$, we say that $N$ has \textbf{w-spectral gap} in $M$ if $N'\cap M^\u=(N'\cap M)^\u$.
\end{defn}

It is interesting to note that a recent result of Tan \cite{tan} shows that a separable II$_1$ factor has property (T) if and only if it has w-spectral gap in every extension, answering a question raised by the first author in \cite{fcep}.  (The original definition of property (T) is more complicated, using the notion of a bimodule over a II$_1$ factor.) 

w-spectral gap is a definability property, as was made explicit by the first author in \cite{spgap}.  Note indeed that we can reformulate the above definition by saying that $N$ has w-spectral gap in $M$ if and only if:  for every $\epsilon>0$, there is a finite subset $F\subseteq N$ and $\delta>0$ such that, for all $x\in B_1(M)$, if $\|xy-yx\|_2<\delta$ for all $y\in F$, then there is $x'\in N'\cap M$ such that $\|x-x'\|_2\leq \epsilon$.  By replacing the condition ``for all $x\in B_1(M)$'' with the stronger condition ``for all $x\in M$'' in the previous sentence, we arrive at the definition of a \textbf{spectral gap} subalgebra; the reason for the terminology can also be found in \cite{spgap}.

In \cite{spgap}, the following fact was established about w-spectral gap subalgebras of e.c. II$_1$ factors:

\begin{thm}\label{bicomm}
Suppose that $N$ is a w-spectral gap subalgebra of the e.c. II$_1$ factor $M$.  Then $(N'\cap M)'\cap M=N$.
\end{thm}

As alluded to in the the previous section, Theorem \ref{bicomm} can be used to give an alternate proof that there is no model companion for the theory of tracial von Neumann algebras.

Using the previous ``bicommutant'' theorem, it is easily verified that if $N$ is a w-spectral gap sub\emph{factor} of the e.c. II$_1$ factor $M$, then $N'\cap M$ is itself a factor, whence so is $N'\cap M^\u$ (as it coincides with $(N'\cap M)^\u$).  Since every II$_1$ factor embeds into some e.c. factor and e.c. factors are locally universal (see the introduction to the first author's article in this volume for a proof of this latter fact), we have arrived at a weak version of Theorem \ref{goldfcep} above, namely every property (T) factor embeds in some locally universal II$_1$ factor with a factorial relative commutant.

In order to complete the proof of Theorem \ref{goldfcep} above, we need to find a common locally universal II$_1$ factor which works for all separable property (T) factors.  \emph{If} Question \ref{eceequestion} above has a positive solution, that is, if all e.c. II$_1$ factors are elementarily equivalent, then we would be finished with the proof. Indeed, under this assumption, any e.c. II$_1$ factor $M$ would satisfy the conclusion of Theorem \ref{goldfcep}.  To see this, fix a property (T) factor $N$ and embed it in an e.c. factor $Q$.  By assumption, $Q^\u\cong M^\u$, whence composing with any isomorphism $Q^\u\to M^\u$, the embedding $N\subseteq Q\hookrightarrow Q^\u$ yields an embedding $N\hookrightarrow M^\u$ with factorial relative commutant. (We have assumed the continuum hypothesis here to ensure that $Q^\u\cong M^\u$; it would be interesting to see if this dependence on the continuum hypothesis can be removed.)

Whether or not Question \ref{eceequestion} has a positive solution, there is a certain subclass of the e.c. factors, namely the class of \textbf{infinitely generic} factors (see Chifan, Drimbe, and Ioana's article in this volume for a definition), which are all pairwise elementarily equivalent.  The class of infinitely generic structures was introduced by Abraham Robinson in \cite{robinson} and imported to the continuous setting in \cite[Section 5]{ecfactors}.  Since every separable tracial von Neumann algebra embeds in a separable infinitely generic factor, running the above argument, but replacing e.c. factors by infinitely generic factors, yields the desired proof of Theorem \ref{goldfcep}.

We point out that recent work of Chifan, Drimbe, and Ioana \cite[Theorem 6.2]{CDI} extended Theorem \ref{goldfcep} by showing that, for any infinitely generic II$_1$ factor $M$ and any separable II$_1$ factor $N$ \emph{without property Gamma}, there is an embedding $i:N\hookrightarrow M^\u$ for which $i(N)'\cap M^\u$ is a factor.

The proof of Theorem \ref{goldfcep} given above points to how one might resolve Conjecture \ref{popafcep}.  Indeed, if one can solve two particular problems in the affirmative, then one can indeed accomplish this goal.

The first problem is purely model-theoretic:

\begin{question}\label{prob1}
Is $\R$ an infinitely generic model of $\Th_\forall(\R)$?
\end{question}

A positive answer to this question was claimed in \cite[Proposition 5.21]{ecfactors}, but the proof there is incorrect.  We believe this question is of independent interest, for a negative solution implies a negative solution to Question \ref{eceequestion} above as it would show that $\R$ is not elementarily equivalent to any infinitely generic model of $\Th_\forall(\R)$.  Indeed, if $M$ is an infinitely generic model of $\Th_\forall(\R)$ and $\R\equiv M$, then since $\R$ is a prime model of its theory (as any embedding of it into its ultrapower is elementary, as discussed above), we have that $\R$ itself is an infinitely generic model of $\Th_\forall(R)$ (as elementary substructures of infinitely generic structures are themselves infinitely generic \cite[Proposition 5.17]{ecfactors}). 

The second problem is purely operator-algebraic and has already been mentioned above in Question \ref{amalgamated}: if $M$ and $N$ are both models of $\Th_\forall(\R)$ with a common subalgebra $Q$, is the amalgamated free product $M*_QN$ also a model of $\Th_\forall(\R)$?

How does this question arise in connection with adapting our proof of Theorem \ref{goldfcep} to settle Conjecture \ref{popafcep}?  Well, the proof of Theorem \ref{bicomm} uses amalgamated free products in a seemingly essential way.  More precisely, if one assumes, towards a contradiction, that there is $b\in ((N'\cap M)'\cap M)\setminus N$, then one obtains a contradiction by considering the embedding $M\subseteq M*_N (N\otimes L(\mathbb{Z}))\hookrightarrow M^\u$ one gets from the fact that $M$ is e.c. and by considering the fact that $b$ commutes with the image of the canonical unitary generator of $L(\mathbb{Z})$ under this embedding.  To run this argument in the relative setting of models of $\Th_\forall(\R)$, one could only quote the fact that $M$ is an e.c. model of $\Th_\forall(\R)$ if the amalgamated free product is itself a model of $\Th_\forall(\R)$.  Note that this argument shows that it suffices to settle this amalgamated free product question in the affirmative when the algebra being amalgamated over has property (T).

In \cite{few}, the authors showed that a positive solution to this amalgamated free product question can be combined with the weaker statement that $\R$ and the infinitely generic models of $\Th_\forall(R)$ have the same \emph{three}-quantifier theory in order to obtain a positive solution to Conjecture \ref{popafcep}.

In a similar spirit to the FCEP but in a somewhat different direction, in Jekel's article in this volume, a free probability-theoretic criterion for a separable tracial von Neumann algebra $M$ to admit an embedding into a matrix ultraproduct $\prod_\u M_n(\bb C)$ with trivial (and thus factorial) relative commutant is established.

\subsection{II$_1$ factors with the generalized Jung property}

The starting point for the main result of this section is Jung's theorem (Fact \ref{jungfact} from above).  Jung's theorem asserts that whenever $N$ is a separable II$_1$ factor that is a model of $\Th_\forall(\R)$ and for which all of its embeddings into $\R^\u$ are unitarily conjugate, then $N$ must be isomorphic to $\R$.  One might also ask what happens instead if one assumes that all embeddings of $N$ into \emph{its own ultrapower} $N^\u$ are unitarily conjugate?  It turns out that, for models of $\Th_\forall(\R)$, this once again characterizes $\R$, as shown by Atkinson and Kunnawalkam Elayavalli in \cite[Corollary 2.7]{AK}: 

\begin{thm}[Atkinson and Kunnawalkam Elayavalli]\label{AK}
If $N$ is a separable II$_1$ factor that is a model of $\Th_\forall(\R)$ for which any two embeddings $N\hookrightarrow N^\u$ are unitarily conjugate, then $N\cong \R$.
\end{thm}

The idea of the proof of the previous result is as follows.  Fix an embedding $\sigma:N\hookrightarrow \R^\u$ and view $\sigma$ as an embedding $\sigma:N\hookrightarrow N^\u$ by composing it with the natural inclusion $\R^\u\subseteq N^\u$.  By assumption, there is a unitary $u=(u_n)_\u\in U(N^\u)$ conjugating $\sigma$ to the diagonal embedding.  Let $E_\R:N\to \R$ denote the canonical conditional expectation map.  Then the map $x\mapsto E_{\R}(u_nxu_n^*):N\to \ell^\infty(\R)$ is a ucp lift of $\sigma$.  Since $\ell^\infty(\R)$ is injective, it follows that $\sigma(N)$ and hence $N$ is injective.  By Connes' fundamental theorem \cite{connes}, this implies that $N\cong \R$.

Using model-theoretic techniques, Atkinson, Kunnawalkam Elayavalli and the first author were able to generalize the previous theorem by merely assuming that any two embeddings of $N$ into its ultrapower $N^\u$ were conjugate by some automorphism of $N^\u$.

\begin{thm}[Atkinson, Goldbring, and Kunnawalkam Elayavalli]\label{AGKtheorem}
If $N$ is a separable II$_1$ factor that is a model of $\Th_\forall(\R)$ for which any two embeddings $N\hookrightarrow N^\u$ are conjugate \emph{by some (not necessarily inner)} automorphism, then $N\cong \R$.
\end{thm}

The proof proceeds in two steps.  The first step is to show that any $N$ satisfying the hypotheses of the theorem must be elementarily equivalent to $\R$.  To see this, first note that the assumption on $N$ implies that any embedding $N\hookrightarrow N^\u$ is elementary. Now consider embeddings $i:\R\hookrightarrow N$ and $j:N\hookrightarrow \R^\u$.  Note then that $j\circ i:\R\hookrightarrow \R^\u$ and $i^\u\circ j:N\hookrightarrow N^\u$ are elementary as are $j^\u\circ i^\u:\R^\u\to (\R^\u)^\u$ and $(i^\u)^\u\circ j^\u:N^\u\hookrightarrow (N^\u)^\u$.  Continuing in this way, we obtain a sequence
$$\R\hookrightarrow N\hookrightarrow \R^\u\hookrightarrow N^\u\hookrightarrow (\R^\u)^\u\hookrightarrow \cdots$$ such that the composition of any two embeddings is an elementary map $\R^{k\u}\hookrightarrow\R^{ (k+1)\u}$ or $N^{k\u}\hookrightarrow N^{(k+1)\u}$, where $\R^{k\u}$ denotes the $k^{\text{th}}$ iterated ultrapower of $\R$ and similarly for $N^{k\u}$.  Letting $N_\infty$ denote the limit of this chain, we see that $\R\preceq N_\infty$ and $N\preceq N_\infty$, whence $\R\equiv N$, as desired.

The second step is to establish the following fact, which is of independent interest:

\begin{thm}\label{eefcep}
If $N\equiv \R$ is separable, then $N$ satisfies the conclusion of Popa's FCEP, that is, there is an embedding $i:N\hookrightarrow \R^\u$ such that $i(N)'\cap \R^\u$ is a factor.
\end{thm}

Before proving Theorem \ref{eefcep}, we make three remarks.  First, this theorem, combined with the first part of the proof, does indeed establish Theorem \ref{AGKtheorem}.  Indeed, if $N$ satisfies the hypotheses of Theorem \ref{AGKtheorem}, then $N\equiv \R$, whence admits an embedding $N\hookrightarrow \R^\u$ with factorial relative commutant.  However, since any two such embeddings are conjugate by an automorphism, we see that all embeddings of $N\hookrightarrow \R^\u$ have factorial relative commutant.  As mentioned in the previous subsection, the embeddings of $N$ into $\R^\u$ with factorial relative commutant are the extreme points of Brown's space $\bb Hom(N,\R^\u)$, whence, under the present assumptions, every point is extreme.  Consequently, there is a single element of $\bb Hom(N,\R^\u)$, whence $N\cong \R$ by Fact \ref{jungfact}.

The second remark is that Theorem \ref{eefcep}, together with Theorem \ref{manymodels} above, gives continuum many nonisomorphic separable examples of II$_1$ factors satisfying Popa's FCEP.

The final remark is that the authors showed in \cite{few}, using Ehrenfeucht-Fra\"isse games, that if $N$ is a II$_1$ factor with the same \emph{four-quantifier theory} as $\R$, then $N$ satisfies the conclusion of Popa's FCEP.  Since we conjecture that $\R$ does not satisfy any sort of quantifier simplification, this should lead to further examples.

We now turn to the proof of Theorem \ref{eefcep}.  The main interest in the proof of this result is that it appears to be the first time that the model-theoretic notion of heir has been used in applications of continuous model theory.

\begin{defn}
Suppose that $T$ is an $L$-theory, $M\models T$, and $A\subseteq B\subseteq M$ are parameter sets.  We say that $q\in S(B)$ is an \textbf{heir} of $p\in S(A)$ if $p\subseteq q$ and for all $L$-formulae $\varphi(x,y)$, parameters $b\in B$, and $\epsilon>0$, there are parameters $a\in A$ such that $|\varphi(x,a)^p-\varphi(x,b)^q|<\epsilon$.
\end{defn}

In the above definition, recall that $\varphi(x,a)^p$ is the value the type $p$ assigns to the formula $\varphi(x,a)$ and likewise for $\varphi(x,b)^q$.

A slightly more complicated version of the usual existence of heir argument can be used to show the following (which was first claimed by Ben Yaacov in \cite{ben} but proven in detail in \cite[Fact 2.2.6]{AGK}):

\begin{lem}
Suppose that $\bb M\models T$ is $\aleph_1$-saturated, $N\preceq \bb M$ is a separable elementary substructure, and $N\subseteq B\subseteq \bb M$ is a separable parameter set.  Then every $p\in S(N)$ has an extension to an element of $S(B)$ that is an heir of $p$. 
\end{lem}

Continuing the proof of Theorem \ref{eefcep}, suppose that $N\subseteq \R^\u$ is an arbitrary separable subalgebra.  By yet another theorem of Brown \cite[Theorem 6.9]{nate}, there is a separable subfactor $N\subseteq M\subseteq \R^\u$ such that $M'\cap \R^\u$ is a factor.  Ordinarily, this would not imply that $N'\cap \R^\u$ is a factor.  However, if $N$ is an \emph{elementary} substructure of $\R^\u$, then one can indeed conclude that $N'\cap \R^\u$ is a factor from the fact that $M'\cap \R^\u$ is a factor.  The key observation is the following:

\begin{prop}
Suppose that $N\subseteq M\subseteq \R^\u$ are separable subfactors and consider types $p\in S(N)$ and $q\in S(M)$ with $q$ an heir of $p$.  If $p(\R^\u)\subseteq N'\cap \R^\u$, then $q(\R^\u)\subseteq M'\cap \R^\u$. 
\end{prop}

The proof of the proposition is routine:  arguing by contrapositive, if there is $b\in M$ such that $\|[x,b]\|_2^q\geq \epsilon$, then there is $a\in N$ such that $\|[x,a]\|_2^p\geq \frac{\epsilon}{2}$.  

It remains to note how the previous proposition establishes Theorem \ref{eefcep} above.  Indeed, we may assume that $N\preceq \R^\u$ and fix a separable subfactor $N\subseteq M\subseteq \R^\u$ such that $M'\cap \R^\u$ is a factor.  We wish to show that $N'\cap \R^\u$ is a factor.  Take $a\in Z(N'\cap \R^\u)$; we  wish to show that $a\in \bb C$.  Let $p:=\tp(a/N)$ and let $q\in S(M)$ be an heir of $p$ to $M$.  Then $q(\R^\u)\subseteq M'\cap \R^\u$ by the previous proposition and $q(\R^\u)\subseteq p(\R^\u)\subseteq Z(N'\cap \R^\u)$ (the latter inclusion follows by noting that any two realizations of $p$ are conjugate by an automorphism of $\R^\u$ fixing $N$ pointwise).  Consequently, $q(\R^\u)\subseteq Z(M'\cap \R^\u)=\bb C$, whence $d(x,\lambda)^q=0$ for some $\lambda\in \bb C$.  It follows that $a=\lambda$, as desired.

We remark that Theorem \ref{eefcep} above yields yet another alternate proof (assuming the continuum hypothesis) that $\Th_\forall(\R)$ has no model companion.  Indeed, suppose towards a contradiction that $\Th(\R)$ is the model companion of its universal theory.  Take a nonstandard model $N$ of $\Th(\R)$.  By assumption, any embedding $N\hookrightarrow N^\u$ is elementary. A standard back-and-forth argument then shows that all embeddings $N\hookrightarrow N^\u$ are conjugate by an automorphism, contradicting Theorem \ref{eefcep}.

By the negative solution to the CEP, there are II$_1$ factors $N$ that are counterexamples to CEP such that any two embeddings $N\hookrightarrow N^\u$ are conjugate by an automorphism.  Indeed, any finitely generic II$_1$ factor will have this property (see \cite{enforceable}).  However, the unitary conjugacy version of this fact (that is, the generalization of Theorem \ref{AK} above to the unrestricted situation) is still open:

\begin{question}\label{Runique?}
Is there a separable II$_1$ factor $N$ that is not isomorphic to $\R$ with the property that any two embeddings $N\hookrightarrow N^\u$ are unitarily conjugate?
\end{question}

Recently, it was shown by Chifan, Drimbe, and Ioana \cite[Theorem 6.4]{CDI} that if $M$ is an infinitely generic factor, then for any \emph{elementary} embedding $i:M\hookrightarrow M^\u$, one has that $i(M)'\cap M^\u$ is a factor.  If there was an infinitely generic factor $M$ for which all embeddings $M\hookrightarrow M^\u$ are elementary, then by the result mentioned in the previous sentence and \cite[Proposition 5.3]{atkinson} (a generalization of a result of Brown referred to above), one would have that any two embeddings of $M$ into $M^\u$ are unitarily conjugate, yielding a positive answer to Question \ref{Runique?}.  Recalling that all embeddings of a finitely generic factor into its ultrapower are elementary, a positive answer to the following question (which is the unrestricted version of Question \ref{prob1}) would consequently lead to a positive answer to Question \ref{Runique?}:

\begin{question}
Do the finitely generic and infinitely generic factors have the same theory?  In other words, is every finitely generic factor also infinitely generic?
\end{question}

While we are on the topic of finitely generic factors, since Theorem \ref{AGKtheorem} above implies that $\R$ is the unique finitely generic separable embeddable factor, it makes sense to ask if this is true in the unrestricted case:

\begin{question}
Is there a unique finitely generic separable factor?
\end{question}

The following general automorphism version of Jung's theorem is also still open:

\begin{question}
Suppose that $N$ is a separable II$_1$ factor that is a model of $\Th_\forall(\R)$ such that any two embeddings $N\hookrightarrow \R^\u$ are conjugate by some (not necessarily inner) automorphism of $\R^\u$.  Must we have $N\cong \R$?
\end{question}

Note that, in model-theoretic terms, the assumption on $N$ in the previous question is simply that the quantifier-free type of $N$ in $\R^\u$ determines its complete type.

\end{document}